\documentstyle[amssymb,12pt]{amsart}

\ifx\mathrm\undefined\let\mathrm\bf\fi
\ifx\mathbf\undefined\let\mathbf\bf\fi

\let\leq\leqslant

\let\al\alpha

\let\gm\gamma 
\let\dl\delta \let\Dl\Delta
 \let\eps\varepsilon \let\epsilon\eps

\let\la\lambda 
\let\om\omega \let\Om\Omega
 \let\phi\varphi

\let\si\sigma

\newcommand{\R}{{\Bbb R}}
\newcommand{\C}{{\Bbb C}}
   
\newcommand{\Ref}[1]{{$($\ref{#1}$)$}}
\newcommand{\bean}{\begin{eqnarray}}
\newcommand{\eean}{\end{eqnarray}}
\newcommand{\be}{\begin{displaymath}}
\newcommand{\ee}{\end{displaymath}}
\newcommand{\bea}{\begin{eqnarray*}}   
\newcommand{\eea}{\end{eqnarray*}}
\newcommand{\g}{{{\frak g}\,}}

\newcommand{\T}{\otimes}

\newcommand{\vs}{\vspace{1.5\baselineskip}}

\newenvironment{proof}{\noindent{\it Proof\/}:\rm}{$\;\Box$
\par\vs}

\newenvironment{example}
{\noindent{\bf Example.\/}}{$\;\Box$ \par\vs}
\newtheorem
{thm}{Theorem}

\newtheorem
{conjecture}{Conjecture}

\newcommand{\End}{{\operatorname{End}}}

\begin{document}
\title{Remarks on critical points of phase functions and norms of Bethe vectors}
\author{Evgeny Mukhin and Alexander Varchenko}
\maketitle

\centerline{\it Dedicated to Peter Orlik on his sixtieth birthday}
\bigskip

\centerline{October 1998}

\begin{abstract}
We consider a tensor product of a Verma module and the linear representation of $sl(n+1)$. We prove that the corresponding phase function, which is used in the solutions of the KZ equation with values in the tensor product, has a unique critical point and show that the Hessian of the logarithm of the phase function at this critical point equals the Shapovalov norm of the corresponding Bethe vector. 
\end{abstract}

\section{Introduction}
Let $\g$ be a simple Lie algebra with simple roots $\al_i$ and Chevalley generators $e_i, f_i, h_i$, $i=1,\dots,n$.
Let $V_1,V_2$ be representations of $\g$  
with highest weights $\la_1,\la_2$. The Knizhnik-Zamolodchikov (KZ) 
equation on a function $u$ with values in $V_1\T V_2$ has the form
\be
\kappa\;\frac{\partial}{\partial z_1} u=\frac{\Om}{z_1-z_2}u, \qquad
\kappa\;\frac{\partial}{\partial z_2} u=\frac{\Om}{z_2-z_1}u,
\ee
where $\Om\in\End(V_1\T V_2)$ is the Casimir operator. Solutions with values in the space 
of singular vectors of weight $\la_1+\la_2-\sum_{i=1}^n l_j\alpha_j$ are given by 
hypergeometric integrals with $l=\sum_{i=1}^n l_j$ integrations, see \cite{SV}.

For an ordered set of numbers $I=\{i_1,\dots,i_m\}$, $i_k\in\{1,\dots,n\}$, 
and a vector $v$ in a representation of $\g$,
denote $f^Iv= f_{i_1}\dots f_{i_m}v$.
The hypergeometric solutions of the KZ equation have the form 
\be
u=\sum u_{I,J} f^Iv_1\T f^Jv_2,\qquad u_{I,J}=\int_{\gm}\Om\tilde{\om}_{I,J}dt_1\wedge\dots\wedge dt_l,
\ee
where $v_1,v_2$ are highest weight vectors of $V_1,V_2$;
the summation is over all pairs of ordered sets $I,J$, such that their union 
$\{i_k, j_s\}$ contains a number $i$ exactly $l_i$ times, $i=1,\dots,n$; $\gm$ is a suitable
cycle; $\tilde{\om}_{I,J}=\tilde{\om}_{I,J}(z_1,z_2,t_1,...,t_l)$ 
are suitable rational functions, the function $\Om=\Om(z_1,z_2,t_1,...,t_l)$,
called the phase function, is given by
\be
\Om = (z_1-z_2)^{(\la_1,\la_2)/\kappa}\;\prod_{j=1}^l (t_j-z_1)^{-(\la_1,\alpha_{t_j})/\kappa}
(t_j-z_2)^{-(\la_2,\alpha_{t_j})/\kappa}
\prod_{1\leq i<j\leq n}(t_i-t_j)^{(\alpha_{t_i},\alpha_{t_j})/\kappa}.
\ee
Here $(\;,\;)$ is the Killing form and $\alpha_{t_i}$ denotes the simple root assigned a
 the variable $t_i$ by the following rule. The first $l_1$ variables $t_1,\dots,t_{l_1}$ are assigned to the simple root $\alpha_1$, the next $l_2$ variables $t_{l_1+1},\dots,t_{l_1+l_2}$ to the second simple root $\alpha_2$, and so on.

Define the normalized phase function $\Phi$ by the formula
\bean\label{Phi}
\Phi(\la_1,\la_2,\kappa)=\prod_{j=1}^l t_j^{-(\la_1,\alpha_{t_j})/\kappa}
(1-t_j)^{-(\la_2,\alpha_{t_j})/\kappa}
\prod_{1\leq i<j\leq n}(t_i-t_j)^{(\alpha_{t_i},\alpha_{t_j})/\kappa}.
\eean
We also substitute $z_1=0$, $z_2=1$ in the rational functions $\tilde{\om}_{I,J}$ and denote the result $\om_{I,J}$.

\begin{conjecture}\label{computable}
If the space of singular vectors of weight $\la_1+\la_2-\sum_{i=1}^n l_j\alpha_j$ is one-dimensional, 
then there is a region $\Dl$ of the form $\Dl=\{t\in\R^l\;|\; 0<t_{\si_l}<\dots<t_{\si_1}<1\}$ 
for some permutation $\si$, such that the integral  $\int_{\Dl}\Phi dt$
can be computed explicitly 
and it is equal to an alternating product of Euler $\Gamma$-functions up to a rational number independent on $\la_1,\la_2,\kappa$.
\end{conjecture}

\begin{example} {\it The Selberg integral.}
Let $\g=sl(2)$. Let $V_1$ and $V_2$ be  $sl(2)$ modules with highest weights 
$\la_1, \la_2 \in \C$. Then the normalized phase function \Ref{Phi} has the form
\bean\label{Selberg phase}
\Phi(\la_1,\la_2,\kappa)=\prod_{j=1}^l t_j^{-\la_1/\kappa} (1-t_j)^{-\la_2/\kappa}\prod_{1\leq i<j\leq l}(t_i-t_j)^{2/\kappa}.
\eean
Conjecture \ref{computable} holds for $\g=sl(2)$ according to the Selberg formula
\be
l!\int_\Dl\Phi(\la_1,\la_2,\kappa)dt_1\dots dt_l=
\prod_{j=0}^{l-1}\,\frac{\Gamma((-\la_1+j)/\kappa+1)\Gamma((-\la_2+j)/\kappa+1)\Gamma((j+1)/\kappa+1)}{\Gamma((-\la_1-\la_2+(2l-j-2))/\kappa+2)\Gamma(1/\kappa+1)},
\ee
where $\Dl=\{t\in\R^l\;|\; 0<t_1<\dots<t_l<1\}$.
\end{example}

Using the phase function $\Phi$ and the rational functions $\om_{I,J}$, one can construct singular vectors in $V_1\T V_2$. Namely, if $t^0$ is a critical point of the function $\Phi$, then the vector $\sum\om_{I,J}(t^0)f^Iv_1\T f^Jv_2$ is singular, see \cite{RV}. The equation for critical points, 
$d\Phi=0$, is called 
the {\it Bethe equation} and the corresponding singular vectors are called the {\it Bethe vectors}.

\begin{conjecture}\label{! critical}
If the space of singular vectors of a given weight in $V_1\otimes V_2$
is one-dimensional, then the corresponding phase function has exactly one critical point modulo permutations of variables assigned to the same simple root.
\end{conjecture}

\begin{example}
\rm
The conjecture holds for $\g=sl(2)$. If $(t_1,\dots,t_l)$ is a critical point  of the function $\Phi(\la_1,\la_2,\kappa)$ given by \Ref{Selberg phase}, then 
\be
\si_k(t)={l \choose k}\prod_{j=1}^k\frac{\la_1-l+j}{\la_1+\la_2-2l+j+1},
\ee
where $\si_1(t)=\sum t_j$, $\si_2(t)=\sum t_it_j$, etc, are the standard symmetric functions, see \cite{V}, so there is a unique critical point up to permutations of coordinates.
\end{example}

The rational functions $\om_{I,J}(t)$ are invariant with respect to permutation of variables assigned to the same simple root. Thus, Conjecture \ref{! critical} implies that there is a unique Bethe vector $X$.

The space $V_1\T V_2$ has a natural bilinear form $B$, called the Shapovalov form, which is the tensor product of Shapovalov forms of factors.
\begin{conjecture}\label{length=hess}
The length of a Bethe vector $X$ equals the Hessian of the logarithm of the phase function $\Phi$ with $\kappa=1$ at a critical point $t^0$,
\be
B(X,X)=det\left(\frac{\partial^2}{\partial t_i\partial t_j}\ln\Phi(t^0)\right).
\ee
\end{conjecture}

\begin{example}
The conjecture holds for $\g=sl(2)$, see \cite{V}.
\end{example}

In this paper we prove Conjectures \ref{computable}, \ref{! critical} and \ref{length=hess} for 
the case when $\g=sl(n+1)$, $V_1$ is a Verma module and $V_2$ is the linear representation.

\section{The integral}

Let
\bean\label{Phi tilde} 
\tilde{\Phi}_n(\alpha,\beta)=t_1^{\alpha_1}(1-t_1)^{\beta_1}\;
\prod_{j=2}^n t_j^{\alpha_j}(t_j-t_{j-1})^{\beta_j}.
\eean

\begin{thm}\label{integral}
Let $\alpha_i>0,\;\beta_i>0$, $i=1,\dots n$. Then
\be
\int_{\Dl_n}\tilde{\Phi}_n(\alpha,\beta)dt_1\dots dt_n=\prod_{j=1}^n\frac{\Gamma(\beta_j+1)\Gamma(\alpha_j+\dots+\alpha_n+\beta_{j+1}+\dots+\beta_n+n-j+1)}{\Gamma(\alpha_j+\dots+\alpha_n+\beta_{j}+\dots+\beta_n+n-j+2)},
\ee
where $\Dl_n=\{t\in\R^n\;|\; 0<t_n<\dots<t_1<1\}$.
\end{thm}
\begin{proof}
The formula is clearly true for $n=1$.

Fix $t_1,\dots,t_{n-1}$ and integrate with respect to $t_n$. We obtain the recurrent relation
\bea
\lefteqn{\int_{\Dl_n}\tilde{\Phi}_n(\alpha,\beta)dt_1\dots dt_n=
\frac{\Gamma(\alpha_n+1)\Gamma(\beta_n+1)}{\Gamma(\alpha_n+\beta_n+2)}\times}\\
&&
\times\int_{\Dl_{n-1}}\tilde{\Phi}_{n-1}(\alpha_1,\dots,\alpha_{n-1},\beta_1,\dots,\beta_{n-2},\beta_{n-1}+\beta_n+\alpha_n+1)dt_1\dots dt_{n-1},
\eea
which implies the Theorem.
\end{proof}

\section{The critical point}

Let $\g=sl(n+1)$. Let $V_1$ be a Verma module of highest weight $\la$, $(\la,\alpha_i)=\la_i$. Let $V_2$ be the linear representation, that is the irreducible representation with highest weight $\om$, $(\om,\alpha_i)=\dl_{i,1}$. 

The nontrivial subspaces of singular vectors of a given weight in the tensor product $V_1\T V_2$ are one dimensional and have weights $\la+\om-\sum_{i=1}^k\alpha_i$, $k=0,\dots,n$.  The computations for weights $\la+\om-\sum_{i=1}^k\alpha_i$, $k<n$, are reduced to the case $\g=sl(k+1)$. Consider the normalized phase function $\Phi_n(\la,\kappa)$ corresponding to the 
weight $\la+\om-\sum_{i=1}^n\alpha_i$. 

We have $\Phi_n(\la,\kappa)=\Phi(\la,\om,\kappa)$, where $\Phi(\la,\om,\kappa)$ is given by \Ref{Phi}. Note that
\be
\Phi_n(\la,\kappa)=\tilde{\Phi}_n(-\la_1/\kappa,\dots,-\la_n/\kappa,-1/\kappa,\dots,-1/\kappa),
\ee
where $\tilde{\Phi}_n$ is given by \Ref{Phi tilde}.

\begin{thm}\label{critical}
The function $\Phi_n(\la,\kappa)$ has exactly one critical point $t^n=(t_1^n,\dots,t^n_n)$ given by 
\be
t_j^n(\la_1,\dots,\la_n)=\prod_{i=1}^j\frac{\la_i+\dots+\la_n+n-i}{\la_i+\dots+\la_n+n-i+1}.
\ee
\end{thm}
\begin{proof}
The computation is obvious if $n=1$.

The equation ${\partial\Phi_n}/{\partial t_n}=0$ has the form
\be
t^n_n=\frac{\la_n}{\la_n+1}t^n_{n-1}.
\ee
Substituting for $t_n^n$ in the equations  ${\partial\Phi_n}/{\partial t_i}=0$, $i=1,\dots,n-1$ and comparing the result with the equation $d\Phi_{n-1}=0$, we obtain
\be
t^n_k(\la_1,\dots,\la_n)=t_k^{n-1}(\la_1\dots,\la_{n-2},\la_{n-1}+\la_n+1),\qquad k=1,\dots,n-1.
\ee
This recurrent relation implies the Theorem.
\end{proof}

\section{The norm of the Bethe vector}

Let $V$ be a $\g$ module with highest weight vector $v$. The Shapovalov form $B(\;,\;): V\T V\to \C$ is the unique symmetric bilinear form with the properties
\be
B(e_ix,y)=B(x,f_iy),\qquad B(v,v)=1,  
\ee
for any $x,y\in V$. The Shapovalov form on a tensor product of 
 modules is the tensor product of Shapovalov forms of factors.

Let $\g=sl(n+1)$. Let $V_1=V_\la$ be a Verma module of highest weight $\la$. Let $V_2=V_\om$ be the linear representation. Then the space of singular vecors in $V_\la\T V_\om$ of weight $\la+\om-\sum_{i=1}^n\alpha_i$ is one-dimensional and is spanned by the Bethe vector $X^n(\la)$ corresponding to the critical point of the function $\Phi_n(\la,\kappa)$. The Bethe vector has the form
\be
X^n(\la)=x_0^n\T f_n\dots f_1v_0+x_1^n\T f_{n-1}\dots f_1v_0+\dots+x_n^n\T v_0,
\ee
where $x^n_i\in V_\la$ and $v_0$ is the highest weight vector in $V_\om$.
Here, $x_0^n=a^nv_\la$, where $v_\la$ is the highest weight vector in $V_\la$ and $a^n$ is the value of the corresponding rational function 
\be
\om_{\emptyset,(n,n-1,\dots,1)}(t)=\frac{1}{t_1-1}\prod_{i=1}^{n-1}\frac{1}{t_{i+1}-t_i}
\ee
at the critical point $t_n$ of function $\Phi_n(\la,\kappa)$, given by Theorem \ref{critical}. For a description of all other rational functions whose values at $t^n$ determine
$x_1^n, ... , x^n_n$, see \cite{SV}. 
We have
\be
a^n=(-1)^n\prod_{k=1}^n\frac{(\la_k+\dots+\la_n+n-k+1)^{n-k+1}}{(\la_k+\dots+\la_n+n-k)^{n-k}}.
\ee

\begin{thm}
\bean\label{norm}
B(X^n(\la),X^n(\la))=\prod_{k=1}^n\frac{(\la_k+\dots+\la_n+n-k+1)^{2(n-k)+3}}{(\la_k+\dots+\la_n+n-k)^{2(n-k)+1}}.
\eean
\end{thm}
\begin{proof}
We also claim
\bean\label{x_n^n norm}
B(x_n^n,x_n^n)=\frac{B(X^n(\la),X^n(\la))}{\la_k+\dots+\la_n+n}.
\eean
Formulas \Ref{norm}, \Ref{x_n^n norm} are readily checked for $n=1$.

The vectors $\{v_0, f_1v_0, f_2f_1v_0,\; \dots,\; f_n \dots f_1v_0\}$ form an orthonormal basis of $V_\om$ with respect to its Shapovalov form.
Clearly, we have
\be
B(X^n(\la),X^n(\la))=\left(\frac{a^n(\la)}{a^{n-1}(\la^\prime)}\right)^2
B(X^{n-1}(\la^\prime),X^{n-1}(\la^\prime))+B(x^n_n,x^n_n),
\ee
where $\la^\prime$ is the $sl(n)$ weight, such that $(\la^\prime,\alpha_i)=\la_{i+1}$, $i=1,\dots,n-1$.

The vector $X^n$ is singular. In particular it means that $e_ix_n^n=0$ for $i>1$ and $e_1x_n^n=-x_{n-1}^n$. 
The vector $x_n^n$ has the form
$x_n^n=\sum_\si b^n_\si f_{\si(1)}\dots f_{\si(n)}v^n_\la$,
where the coefficients $b^n_\si$ are the values of the corresponding rational functions at the critical point given by Theorem \ref{critical}.

Let $b^n=b_{\si={\rm id}}^n$. Then we have
\bea
B(x_n^n,x_n^n)=B(x_n^n,b^n\,f_1\dots f_n v^n_\la)=-b^n\,B(x_{n-1}^n,f_2\dots f_n v^n_\la)=
\\
=-b^n\,\frac{a_n}{a_{n-1}}B(x_{n-1}^{n-1},f_1,\dots f_{n-1}v^{n-1}_{\la^\prime})=-\frac{b^n}{b^{n-1}}\;\frac{a_n}{a_{n-1}}B(x_{n-1}^{n-1},x_{n-1}^{n-1}),
\eea
where $x^{n-1}_{n-1}$ is a component of the singular vector in $V_{\la^\prime}\T V_\om$.

The coefficient $b^n$ is the value of the function 
\be
\om_{(n,n-1,\dots,1),\emptyset}(t)=\frac{1}{t_n}\prod_{i=1}^{n-1}\frac{1}{t_i-t_{i+1}}
\ee
at the critical point $t^n$, given by Theorem \ref{critical}.
We have
\be
b^n=(-1)^{n-1}\frac{a_n}{\la_1+\dots+\la_n+n}\prod_{k=1}^n\frac{\la_k+\dots+\la_n+n-k+1}{\la_k+\dots+\la_n+n-k}
\ee

Now, formulas \Ref{norm}, \Ref{x_n^n norm} are proved by induction on $n$.
\end{proof}

\begin{thm}
\be
B(X^n(\la),X^n(\la))=
det\left(\frac{\partial^2}{\partial t_i\partial t_j}\ln\Phi_n(\la,\kappa=1)(t^n)\right),
\ee
where $t^n$ is the critical point of the phase function $\Phi_n(\la,\kappa)$ given by Theorem \ref{critical}.
\end{thm}
\begin{proof}
It is sufficient to prove the Theorem for $\la_i>0,\kappa<0$.  We tend $\kappa$ to zero and compute the asymptotics of the integral $\int_{\Dl_n}{\Phi_n} dt$.

On one hand, the integral is evaluated by Theorem \ref{integral}. We compute the asymptotics using the Stirling formula for $\Gamma$-functions. 

On the other hand, the asymptotics of the same integral can be computed by the method of stationary phase, since the critical point  $t^n$ of the function $\Phi_n$ is non-degenerate by Theorem 1.2.1 in \cite{V}. Then the asymptotics of the integral is 
\be
(2\pi\kappa)^{l/2} \Phi_n(\la,\kappa)(t^n)\left({\rm Hess}(\kappa\ln\Phi_n(\la,\kappa)(t^n)\right)^{-1/2}.
\ee
Note that $\kappa\ln\Phi_n(\la,\kappa)=\ln\Phi_n(\la,1)$, and
\be
\Phi_n(\la,\kappa)(t^n)=\prod_{k=1}^n\frac{(\la_k+\dots+\la_n+n-k+1)^{(\la_k+\dots+\la_n+n-k+1)/\kappa}}{(\la_k+\dots+\la_n+n-k)^{(\la_k+\dots+\la_n+n-k)/\kappa}}.
\ee

Comparing the results we compute the Hessian explicitly and prove the Theorem.
\end{proof}

\bigskip

{\it Mathematical Sciences Research Institute, 1000 Centennial Drive, Berkeley, CA 94720-5070}

{\it Department of Mathematics, University of North Carolina at Chapel Hill, Chapel Hill, NC 27599-3250, USA.}

{\it E-mail addresses:} {\rm mukhin@@msri.org,
av@@math.unc.edu}

\end {document}